\documentclass[letterpaper,journal,9pt]{IEEEtran}

\IEEEoverridecommandlockouts                              

\usepackage{verbatim,color}

\usepackage{amsmath, amssymb, bbm, xspace}
\usepackage{epsfig}
\usepackage{longtable}
\usepackage{color}
\usepackage{mathrsfs}
\usepackage{subfigure}

\usepackage{courier}

\newtheorem{theorem}{Theorem}[section]

\newtheorem{proposition}[theorem]{Proposition}

\newtheorem{corollary}[theorem]{Corollary}

\newtheorem{definition}{Definition}[section]

\def\bkE{{\rm I\kern-.17em E}}
\def\bk1{{\rm 1\kern-.17em l}}
\def\bkD{{\rm I\kern-.17em D}}
\def\bkR{{\rm I\kern-.17em R}}
\def\bkP{{\rm I\kern-.17em P}}
\def\Cbb{{\mathbb{C}}}

\def\fnat{{\bf F}^{\textrm{nat}}}

\def\bkZ{{\bf{Z}}}

\def\bkE{{\rm I\kern-.17em E}}
\def\bk1{{\rm 1\kern-.17em l}}
\def\bkD{{\rm I\kern-.17em D}}
\def\bkR{{\rm I\kern-.17em R}}
\def\bkP{{\rm I\kern-.17em P}}

\def\bkZ{{\bf{Z}}}
\def\b12{(\beta_1,\beta_2)}
\def\Epec{{\mathscr{E}}}

\newcounter{example}
\renewcommand{\theexample}{\thesection.\arabic{example}}
\newenvironment{examplec}[1][]{\refstepcounter{example}
\par\medskip \noindent%
   \textbf{Example~\theexample. #1} \rmfamily}{\hfill $\square$   \hspace{-4.5pt} \vspace{6pt}}

\newcounter{remark}
\renewcommand{\theremark}{\thesection.\arabic{remark}}
\newenvironment{remarkc}[1][]{\refstepcounter{remark}
\noindent{\itshape Remark~\theremark. #1} \rmfamily}{\hspace*{\fill}~$\square$\vspace{0pt}}

\def\t{^\top}

\newlength{\noteWidth}
\long\def\notes#1{\ifinner
{\tiny #1}
\else
\marginpar{\parbox[t]{\noteWidth}{\raggedright\tiny #1}}
\fi\typeout{#1}}

 \def\notes#1{\typeout{read notes: #1}} 

\def\mi{^{-i}}

\newcommand{\ie}{i.e.\@\xspace} 
\newcommand{\eg}{e.g.\@\xspace} 
\newcommand{\etal}{et al.\@\xspace} 

\newcommand{\Real}{\ensuremath{\mathbbm{R}}}

\newcommand{\dom}{\mathrm{dom}}
\newcommand{\minimize}[1]{\displaystyle\minim_{#1}}
\newcommand{\minim}{\mathop{\hbox{\rm min}}}
\newcommand{\maximize}[1]{\displaystyle\maxim_{#1}}
\newcommand{\maxim}{\mathop{\hbox{\rm max}}}
\DeclareMathOperator*{\st}{subject\;to}

\def\subject{\hbox{\rm s.t.}}

\def\fnat{{\bf F}^{\rm nat}}

\def\Cbb{\mathbb{C}}
\def\Cbb{{\mathbb{C}}}

\def\ind{{\rm ind}}

\def\range{{\rm range}}

\def\half  {{\textstyle{1\over 2}}}

\def\spose#1{\hbox to 0pt{#1\hss}}

\def\text #1{\hbox{\quad#1\quad}}

\def\nthinsp{\mskip -2   mu}

\def\superstar{^{\raise 0.5pt\hbox{$\nthinsp *$}}}
\def\SUPERSTAR{^{\raise 0.5pt\hbox{$*$}}}

\def\lamstarT {\lambda^{\raise 0.5pt\hbox{$\nthinsp *$}T}}

\def\Ascr{{\cal A}}
\def\Bscr{{\cal B}}
\def\Fscr{{\cal F}}
\def\Dscr{{\cal D}}
\def\Iscr{{\cal I}}

\def\Lscr{{\cal L}}

\def\Tscr{{\cal T}}

\def\Pscr{{\cal P}}
\def\Qscr{{\cal Q}}
\def\Sscr{{\cal S}}

\def\Nscr{{\cal N}}
\def\Rscr{{\cal R}}

\def\th{^{\rm th}}

\def\aur{\;\textrm{and}\;}
\def\where{\;\textrm{where}\;}

\def\non{\nonumber}
\def\SOL{{\rm SOL}}
\def\VI{{\rm VI}}

\let\forallnew\forall
\renewcommand{\forall}{\forallnew\ }
\let\forall\forallnew
\newcommand{\mlmfg}{multi-leader multi-follower game\@\xspace}
\newcommand{\mlmfgs}{multi-leader multi-follower games\@\xspace}

		\def\bkE{{\rm I\kern-.17em E}}
		\def\bk1{{\rm 1\kern-.17em l}}
		\def\bkD{{\rm I\kern-.17em D}}
		\def\bkR{{\rm I\kern-.17em R}}
		\def\bkP{{\rm I\kern-.17em P}}
		\def\bkY{{\bf \kern-.17em Y}}
		\def\bkZ{{\bf \kern-.17em Z}}
		\def\bkC{{\bf  \kern-.17em C}}

{\begin{list}{}%
         {\setlength{\leftmargin}{#1}}%
         \item[]%
}
{\end{list}}

		\def\bsp{\begin{split}}
		\def\beq{\begin{eqnarray}}
		\def\bal{\begin{align*}}
		\def\bc{\begin{center}}
		\def\be{\begin{enumerate}}
		\def\bi{\begin{itemize}}
		\def\bs{\begin{small}}
		\def\bS{\begin{slide}}
		\def\ec{\end{center}}
		\def\ee{\end{enumerate}}
		\def\ei{\end{itemize}}
		\def\es{\end{small}}
		\def\eS{\end{slide}}
		\def\eeq{\end{eqnarray}}
		\def\eal{\end{align*}}
		\def\esp{\end{split}}

		\def\problemsmall#1#2#3#4{\fbox
		 {\begin{tabular*}{0.47\textwidth}
			{@{}l@{\extracolsep{\fill}}l@{\extracolsep{6pt}}l@{\extracolsep{\fill}}c@{}}
				#1 & $\minimize{#2}$ & $#3$ & $ $ \\[5pt]
					  $\subject\ $ &    & $#4$ & $ $
			\end{tabular*}}
			}

		\def\problem#1#2#3#4{\fbox
		 {\begin{tabular*}{0.80\textwidth}
			{@{}l@{\extracolsep{\fill}}l@{\extracolsep{6pt}}l@{\extracolsep{\fill}}c@{}}
				#1 & $\minimize{#2}$ & $#3$ & $ $ \\[5pt]
					 & $\subject\ $    & $#4$ & $ $
			\end{tabular*}}
			}

		\def\minmaxproblemcompactsmall#1#2#3#4#5{\fbox
		 {\begin{tabular*}{0.47\textwidth}
			{@{}l@{\extracolsep{\fill}}l@{\extracolsep{6pt}}l@{\extracolsep{\fill}}c@{}}
				#1 & $\minimize{#2}$ $ \maximize{#3}$& $#4$ & 			 $\ \subject$     $#5$ 
			\end{tabular*}}
			}

\def\maxproblem#1#2#3#4{\fbox
		 {\begin{tabular*}{0.80\textwidth}
			{@{}l@{\extracolsep{\fill}}l@{\extracolsep{6pt}}l@{\extracolsep{\fill}}c@{}}
				#1 & $\maximize{#2}$ & $#3$ & $ $ \\[5pt]
					 & $\subject\ $    & $#4$ & $ $
			\end{tabular*}}
			}

	\def\problemcompactsmall#1#2#3#4{\fbox
		 {\begin{tabular*}{0.47\textwidth}
			{@{}l@{\extracolsep{\fill}}l@{\extracolsep{6pt}}l@{\extracolsep{\fill}}c@{}}
				#1 &$\minimize{#2}$\ $#3$ &{\rm s.t.} $#4 $ 
			\end{tabular*}}}
	\def\maxproblemcompactsmall#1#2#3#4{\fbox
		 {\begin{tabular*}{0.47\textwidth}
			{@{}l@{\extracolsep{\fill}}l@{\extracolsep{6pt}}l@{\extracolsep{\fill}}c@{}}
				#1 &$\maximize{#2}$\ $#3$ &{\rm s.t.}\ $#4 $ 
			\end{tabular*}}}

	\def\cp2problem#1#2#3#4{\fbox
		 {\begin{tabular*}{0.9\textwidth}
			{@{}l@{\extracolsep{\fill}}l@{\extracolsep{6pt}}l@{\extracolsep{\fill}}c@{}}
				#1 & & $#4 $ 
			\end{tabular*}}}

		\def\bkE{{\rm I\kern-.17em E}}
		\def\bk1{{\rm 1\kern-.17em l}}
		\def\bkD{{\rm I\kern-.17em D}}
		\def\bkR{{\rm I\kern-.17em R}}
		\def\bkP{{\rm I\kern-.17em P}}
		
		\def\bkZ{{\bf{Z}}}

\newcommand {\beeq}[1]{\begin{equation}\label{#1}}
\newcommand {\eeeq}{\end{equation}}
\newcommand {\bea}{\begin{eqnarray}}
\newcommand {\eea}{\end{eqnarray}}

\def\texitem#1{\par\smallskip\noindent\hangindent 25pt
               \hbox to 25pt {\hss #1 ~}\ignorespaces}

\def\st{\mbox{subject to}}

\usepackage{epsfig,verbatim}
\usepackage{amsmath}
\usepackage{epsfig,amsmath,amssymb}
\usepackage{longtable}
\usepackage{color}
\usepackage{mathrsfs}

\def\texitem#1{\par\smallskip\noindent\hangindent 25pt
               \hbox to 25pt {\hss #1 ~}\ignorespaces}

\def\st{\mbox{subject to}}

\def\PS{{\rm P}^{\rm ae}}

\def\LiS{{\rm L}^{\rm ae}_i}
\def\EpecS{\Epec^{\rm ae}}
\def\FscrS{\Fscr^{\rm ae}}
\def\Iscr{\mathcal I}
\def\bl{\rm bl}

\renewcommand{\emph}[1]{\textbf{#1}}
\def\bkE{{\rm I\kern-.17em E}}
\def\bk1{{\rm 1\kern-.17em l}}
\def\bkD{{\rm I\kern-.17em D}}
\def\bkR{{\rm I\kern-.17em R}}
\def\bkP{{\rm I\kern-.17em P}}

\def\bkZ{{\bf{Z}}}
\def\b12{(\beta_1,\beta_2)}
\def\Epec{{\mathscr{E}}}
\newcommand{\gap}{\vspace{0in}}
\def\us#1{{{\color{black}#1}}}
\def\ae{\rm ae}
\def\cc{\rm cc}
\def\bl{\rm bl}

\def\Pcent{{\rm P}^{\rm quasi}}

\begin{document}
\title{An Existence Result for Hierarchical Stackelberg v/s Stackelberg Games}
\author{Ankur A. Kulkarni and Uday V. Shanbhag\thanks{Ankur is with the
	Systems and Control Engineering group at the Indian Insitute of
		Technology Bombay, Mumbai 400076, India, while Uday is with the
		Department of Industrial and Manufacturing Engineering at the
		Pennsylvania State University at University Park.  They  are
		reachable at ({\tt kulkarni.ankur@iitb.ac.in,udaybag@psu.edu}).
		Ankur's work has been partially funded by the research grant
		associated with the INSPIRE Faculty Award, Department of Science
		and Technology, Government of India. Uday's work has been
		partially funded by the NSF CMMI 124688 (CAREER). 
	}}
\date{}
\maketitle
\begin{abstract}
In Stackelberg v/s Stackelberg games a collection of leaders compete in a Nash game constrained by the equilibrium conditions of another Nash game amongst the followers. The resulting equilibrium problems are plagued by the nonuniqueness of follower equilibria and nonconvexity of leader problems whereby the problem of providing sufficient conditions for existence of global or even local equilibria remains largely open. Indeed available existence statements are restrictive and model specific. In this paper, we present what is possibly the first general existence result for  equilibria 
for this class of games. Importantly, we impose no  single-valuedness
assumption on the equilibrium of the follower-level game.
Specifically, under the assumption that the objectives of the leaders admit
a \textit{quasi-potential} function, a concept we introduce in this paper, the global and local minimizers  of
a suitably defined optimization problem are shown to be the global and
local equilibria of the game.  In effect existence of equilibria can be
guaranteed by {the solvability of} an optimization problem, which {holds
under mild and verifiable conditions}.  We motivate {quasi-potential}
games through an application in communication networks. \end{abstract}

\section{Introduction} \label{sec:intro}
The recent past has seen increased interest in hierarchical systems with
competing participants.  The standard analysis for such systems
concentrates on the Stackelberg model~\cite{vonstack52theory} of a game
-- a single player, called \us{a} \textit{leader}, acts first \us{while}
all other players, called followers, act
	\us{subsequently} under the usual assumptions of a noncooperative game. However
contemporary markets, such as those in the power industry,  require the
modeling and analysis of a more complex game where multiple Stackelberg
leaders compete in a noncooperative game following which followers play
a noncooperative game amongst themselves, taking the decisions of leaders as fixed.  This situation models the
clearing of a sequence of markets, such as the
day-ahead and real-time markets, where the participants of the day-ahead
market are taken as leaders and those of the real-time markets are taken
as followers. 

Let $\Nscr = \{1,2,\hdots,N\}$ denote the set of leaders where leader
$i$ solves the parametrized problem: 
$$
{\small	\problemcompactsmall{L$_i(x^{-i})$}
	{x_i,y_i}
	{\varphi_i(x_i,y_i;x^{-i})}
				 {\begin{array}{r@{\ }c@{\ }l}
		x_i &\in& X_i, 
		y_i \in  \Sscr(x),
	\end{array}}}
	$$
$x_i \in \Real^{m_i}$ and $\varphi_i$ denote leader $i$'s action
and objective while $x^{-i}$ and $(\bar x_i,x^{-i})$ are defined as
$x^{-i} \triangleq
(x_1,\hdots,x_{i-1},x_{i+1},\hdots,x_N) $ and $
(\bar{x}_i,x^{-i}) \triangleq
(x_1,\hdots,x_{i-1},\bar{x}_i,x_{i+1},\hdots,x_N).$ 
 For each $x$,
	the set of follower equilibria \us{is} denoted by $ \Sscr(x)$ and $y_i$
	denotes the strategy profile of \textit{all} followers. Though $y_i$ is not strictly within leader $i$'s control, a minimization over
	$y_i$ is performed by an {\em optimistic} leader; a
	{\em pessimistic} leader would maximize over $y_i$ while minimizing over
	$x_i$. These notions coincide if $ \Sscr $ is single-valued but can be quite different when $ \Sscr $ is multi-valued.
We assume that each follower solves a convex optimization problem parametrized by the strategies of the leaders and other followers. It follows that $\Sscr(x)$ is given by the solution set of a variational inequality (VI), say VI$(G(x,\cdot),K(x))$; the solution set is denoted $\SOL(G(x_i,x^{-i},\cdot),K(x_i,x^{-i}))$.
We assume that the set-valued map $ K $ is continuous and $ G $ is a continous mapping of all variables.

The set $X_i$ represents other constraints and is assumed to
be a convex and compact set.  For each $i$, objective function $\varphi_i$ is defined over $X \times Y$, where $X
\triangleq \prod_{i=1}^N X_i$ and $ Y $ is the ambient space of $ y_i $ for any $ i \in \Nscr $. We assume $ \varphi_i $ to be continuous for each $ i $.  Let $y =
(y_1,\hdots,y_N)$ and   $\Omega_i(x^{-i})$  be the feasible
region of L$_i(x^{-i})$, {given by}
\begin{align}
	\Omega_i(x^{-i}) \triangleq 
	\left\{(x_i,y_i) \ \left\lvert \  x_i  \in X_i, y_i \in \Sscr(x)
\right. \right\}. \label{eq:omegai}\end{align}
Let $\Omega(x)$ denote the Cartesian product of $\Omega_i(x\mi),$ $ \Omega(x) \triangleq \prod_{i=1}^N
\Omega_i(x^{-i}) $ and 
let $\Fscr$ be \us{defined as} \begin{align}
\Fscr \triangleq \left\{(x,y)  \ \left\lvert x_i \in X_i, y_i \in Y_i,
y_i \in \Sscr(x), i \in \Nscr \right. \right\}, \label{eq:fscr} \end{align}
which is the set of tuples $(x,y)$ such that $(x_i,y_i)$ is feasible for L$_i(x\mi)$ for all $i$. 
It is easily seen that $\Fscr$ is the set of fixed points of $\Omega$, \ie $ \Fscr = \{(x,y) \in \Real^n: (x,y) \in \Omega(x,y)\}. $
We denote this multi-leader multi-follower game by $\Epec$. 

In the spirit of the Stackelberg equilibrium, or subgame perfect equilibrium, the followers' response must be taken as the Nash equilibrium of the follower-level game parametrized by leader strategy profiles, whereas the leaders themselves choose decisions their while anticipating and being constrained by the equilibrium of the follower-level game. An equilibrium of such a multi-leader multi-follower game or Stackelberg v/s Stackelberg game is the Nash equilibrium of this game between leaders.

\begin{definition}[Global Nash equilibrium]
	Consider the multi-leader multi-follower game $\Epec$. The {\em global
	Nash equilibrium}, or simply {\em equilibrium}, of $\Epec$ is a point $(x,y) \in \Fscr$ that satisfies the following: 
\begin{align}\label{br}
\varphi_i(x_i,y_i;x^{-i}) &\leq \varphi_i(u_i,v_i;x^{-i}),
\end{align}
$\forall (u_i,v_i) \in \Omega_i(x^{-i},y^{-i}), $ and all $ i \in \Nscr. $
\end{definition}
Eq \eqref{br} says that at an equilibrium $(x,y)$,  $(x_i,y_i)$ lies in the set of best responses to
$(x^{-i},y^{-i})$ for all $i.$ The qualification ``global'' is useful in
distinguishing the equilibrium from its stationary counterparts (referred
		to as a  ``Nash B-stationary'')  or its local
counterpart (referred to as a ``local Nash equilibrium''). In this paper, we mainly  
focus on global Nash equilibria and we refer to  them as simply
``equilibria''; other notions are qualified accordingly.

Despite it being a reasonable and natural model, 
till date no reliable theory for the existence of equilibria to \mlmfgs is available~\cite{pang05quasi}. 
In fact there are fairly
simple \mlmfgs which admit no equilibria; a particularly telling example was shown by Pang and Fukushima~\cite{pang05quasi} (we discuss this example in Section \ref{sec:pot-examp}).  
To the best of our knowledge, existence statements are known only for problems where the follower-level equilibrium  is unique for each strategy profile of the leaders. These results are obtained via explicit substitution of this equilibrium and further analysis of the leader-level equilibrium~\cite{sherali84multiple,su07analysis,demiguel09stochastic}; it is evident that such a modus operandi can succeed only when the cost functions and constraints associated with the players take suitable simple forms. 
On another track, existence has been claimed for weaker notions
of equilibria, \eg,  solutions of  the
aggregated stationarity conditions of the problems of the
leaders~\cite{hu07epec, SIGlynn11,outrata04note}.

The analysis of such games is hindered by the inapplicability of
standard fixed point theorems -- in particular the lack of convexity of
the leaders' problems, and suitable continuity properties in the best
response functions of the leaders.  In this paper we provide a clean
result for the existence of equilibria that does not assume the
uniqueness or favorable structure 
of the follower-level equilibrium. The tool we use is similar to but distinct from the potential game~\cite{monderer96potential}. We relate the minimizer of an optimization problem to the equilibrium of the game, thereby obviating
the need to apply fixed point theory.  
\us{Our main contributions can be summarized as follows:}

\noindent {\bf (1) Global equilibria of quasi-potential games:} {We
introduce a new class of multi-leader multi-follower games, called
\textit{quasi-potential games}, in which the leaders' objectives take on
a particular structure, part of which admits a potential function. We
motivate this class of games through an application in communication
networks. We show that the global minimizers of a suitably defined
optimization problem are the global equilibria of such games.
Consequently, sufficiency requirements for existence reduce to mild and
verifiable conditions for the solvability of optimization problems (such
		as continuity of the objective and compactness of feasible
		region).  Furthermore, we establish a similar relation for the
pessimistic formulation.\\ 
\noindent {\bf (2) Local and Nash-Stationary Equilibria:} Notably, such
relationships are shown to extend to allow for relating local minimizers
and stationary points of the associated optimization problem to local
Nash equilibria and stationary equilibria. }

The remainder of the paper is organized into three sections.  In
Section~\ref{sec:conventional}, we provide an example of
\us{a quasi-potential} \mlmfg and comment on the analytical intractability of general \mlmfgs.
In Section~\ref{sec:standard}, we provide our main results. The paper
concludes in Section~\ref{conc} with a brief summary.

\section{Multi-leader multi-follower games: examples and background} \label{sec:conventional}

\subsection{Examples {of multi-leader multi-follower games}}\label{sec:examples}
In one of the first examples of \mlmfgs~\cite{sherali84multiple}, a set
	of followers compete in a Cournot game to determine quantity of
	production, while each leader \us{makes} a decision constrained by the
	equilibrium quantities produced by the followers. Furthermore,
				leaders compete amongst each other to decide their
					production levels, subject to the levels they
					anticipate from the followers. Following is another
					example from communication networks.
\paragraph{Congestion control in communication networks} Suppose $\Nscr =\{1, \hdots, N\}$ denotes the
set of users and $K_1, \hdots, K_N$ are strategy sets of users. Given a set of
 flow decisions $x = (x_1,\hdots,x_N)$ of users, the network
manager solves a parametrized optimization problem given by the following:
$$ {\problemcompactsmall{Net$(x)$}
	{y}
	{f(y;x) }
				 {\begin{array}{r@{\ }c@{\ }l}
	y & \in & C(x),
\end{array}}}$$
in which $y$ represents the network decisions (which include flow
specifications etc.), $f(y;x)$ denotes the network manager's objective and $C(x)$ represents the set of feasible allocations
available to the manager for user decisions $x$. Each user is assumed to be a leader
with respect to the network manager (follower) and the resulting user
problem is given by the following:
$$ {\maxproblemcompactsmall{L$_i(x^{-i})$\hspace{1mm}}
	{x_i,y}
	{\hspace{1mm}U_i(x_i)-h(y)}
				 {\hspace{-2mm}\begin{array}{r@{\ }c@{\ }l}
					 x_i & \in & K_i, 		y \in  \mbox{SOL(Net}(x)),
\end{array}}}$$
where $h(y)$ represents the congestion cost associated with the network
manager's decision $y$ and $ \SOL(\cdots) $ represents the solution of problem `$ \cdots $'. We assume that every user is charged the entire
cost of congestion, an assumption that is standard in such mdoels
(cf.~\cite{alpcan03distributed,basar07control}). Such a model represents a hierarchical
generalization of the competitive model considered by
Ba\c{s}ar~\cite{basar07control}, Alpcan and
Ba\c{s}ar~\cite{alpcan03distributed}, and Yin, Shanbhag, and
Mehta~\cite{yin09nash2}. Notably, past work has not modeled the network
manager as a separate entity but this extension has much relevance when
considering the role of large, and possibly strategic, independent
\us{service} providers  in the context of network management.
\us{Similar models may arise in power markets where the ``cost of grid
	reliability'' is socialized~\cite{hobbs-strategic}. } 
\subsection{\us{A comment on} the intractability of \mlmfgs} \label{sec:int}
We now briefly comment on analytical difficulties that arise in games such as $ \Epec $. An equilibrium of $ \Epec $ 
is the Nash equilibrium of the game where players solve problems $ \{{\rm L}_i\}_{i\in \Nscr} $. 
At the Nash equilibrium, each player's strategy is his ``best response''
assuming  the strategies of his opponents are held fixed. For any tuple
of strategies $(x,y)$, one may define a {\em reaction map} $\Rscr : \dom(\Omega) \rightarrow 2^{{\footnotesize\range}(\Omega)}$
(set-valued in general) as 
$\Rscr(x,y) := \prod_{i = 1}^N \Rscr_i(x^{-i},y^{-i}), $ where $\Rscr_i(x^{-i},y^{-i}) = \SOL({\rm L}_i(x\mi))$. Although $ \Omega_i(\cdot) $ is independent of $ y\mi $, for this section we will write $ \Omega_i $ explicitly as a function of $ y\mi $ with the identification  $ \Omega_i(x\mi,y\mi) \equiv \Omega_i(x\mi). $
$ (x,y) $ is an equilibrium of $ \Epec $ if and only if $ (x,y) $ is a fixed point of $ \Rscr, $ and one could, in principle, approach the problem through fixed point theory. However due to the nonconvexity of problems $\{ {\rm L}_i\}_{i\in \Nscr}$ 
	difficulties arise when one attempts to apply fixed point theorems to this reaction map.  
Almost all fixed point theorems rely on
the following categories of assumptions: 
\begin{enumerate}
\item  (a) the mapping to which a fixed point is sought is assumed to be a {\em self-mapping}; 
\item (b) (i) the domain of the mapping and (ii) the images are required to be of a specific shape, \eg convex;
\item (c) the mapping is required to be continuous (if the mapping is single-valued) or upper semicontinuous (if set-valued). 
\end{enumerate}
 The first difficulty encountered is that $\Rscr$ is not necessarily a self-mapping: $\Rscr$ maps $\dom(\Omega)$ to 
 $\range(\Omega)$ and $\range(\Omega)$ may not be a subset of
 $\dom(\Omega)$. Second, $\dom(\Omega)$ is hard to characterize and 
little can be said about its shape. Finally, the continuity (or upper
		semicontinuity) of $\Rscr$ is far from immediate.  There are
ways of circumventing difficulties (a) and (b), some of which have been
employed in literature. If $\Omega(x,y) \neq \emptyset$ for all $(x,y)
	\in X \times Y^N$, where $ Y^N = \prod_{i \in \Nscr} Y $, we get $\dom(\Omega) = X \times Y^N$ and $\Rscr$ may
	be taken to be a map from $X \times Y^N$ to subsets of $X \times Y^N$.
	This approach was employed by Arrow and Debreu \cite{arrow54existence}. If $X \times Y^N$ is convex, as is in our case, the
	difficulty (b(i)) is also circumvented. The (upper semi-)continuity
	of $\Rscr$ requires $\Omega$ to be continuous~\cite{facchinei03finiteI}, a property
	that rarely holds if $\Sscr$ (the solution set of a VI) is multivalued.  As a consequence, the
	shape of the mapped values ((b.ii)) and the upper semicontinuity of $\Rscr$
	are the key barriers to the success of this approach. In the case
	where $\Sscr$ is single-valued, the continuity of $\Rscr$ follows
	readily (since $ \Sscr $, the solution set of a parametrized VI is upper-semicontinuous with respect to the parameter). A majority of the known results for \mlmfgs are indeed for 	this case.

\section{Quasi-potential multi-leader multi-follower	games}\label{sec:standard} This section contains our main results. In
	Section~\ref{sec:global_stand}, we define a new class of
	multi-leader multi-follower games, called {\em quasi-potential}
	games,  and provide conditions for the existence of global
	equilibria for the optimistic and pessimistic formulation. Analogous results for local and Nash stationary
	equilibria, for the optimistic formulation are derived in Section~\ref{sec:local_quasi}. Finally,
	in Section~\ref{sec:pot-examp} we conclude by commenting on how the
	results apply to  some well-known examples.

\subsection{Existence results for  global equilibria} \label{sec:global_stand}
To motivate quasi-potential games, let us first consider the case where for
each tuple of leader decisions, there is a unique follower equilibrium.
In this case, the leader problem L$_i$ can be reformulated as follows:
$$
{	\problemcompactsmall{$\tilde{\mbox{L}}_i(x^{-i})$}
	{x_i}
	{\varphi_i(x_i,\Sscr(x);x^{-i})}
				 {\begin{array}{r@{\ }c@{\ }l}
		x_i &\in& X_i.
	\end{array}}}
$$
 Denote this game as $\tilde{\Epec}$.  We may
now define an implicit potential multi-leader multi-follower game as follows. 
\begin{definition}[Implicit Potential multi-leader multi-follower games]\label{def:impot}
An implicit potential multi-leader multi-follower  game is a multi-leader
multi-follower game where for each $ x \in X $, $ |\Sscr(x)|=1$. The $i\th$ leader solves $\tilde{{\rm L}}_i(x^{-i})$ and there exists a function $\pi$ 
	 such that for all $i\in \Nscr$, for all $(x_i,x\mi) \in X$, 
and for all $x_i' \in X_i$,
\begin{align*}
\varphi_i(x_i,\Sscr(x);x\mi)&-\varphi_i(x_i',\Sscr(x_i',x\mi);x\mi)  \\ &=
\pi(x_i,\Sscr(x);x\mi) - \pi(x_i',\Sscr(x_i',x\mi);x\mi). 
\end{align*}
\end{definition}
We thus have the following proposition; the proof is straight-forward, hence we skip it.
\begin{proposition} Consider an implicit potential multi-leader multi-follower
game $\tilde{\Epec}$. Then any global minimizer of P$^{\rm imp}$ is a global Nash
equilibrium of $\tilde{\Epec}$, where P$^{\rm imp}  $ is the following problem:
$$
	\problemcompactsmall{P$^{\rm imp}$}
	{x}
	{\pi(x)}
				 {\begin{array}{r@{\ }c@{\ }l}
					x \in X. 
	\end{array}}
	$$
\end{proposition}

In principle, the above is an existence result. However two challenges
emerge when one considers the result in practice. First, in most cases,
	   particularly when the second-level equilibrium constraints arises
	   from the equilibrium conditions of a constrained problem,
	   $\Sscr(\cdot)$ is not single-valued. Second, even if
	   $\Sscr(\cdot)$ is single-valued, ascertaining the potentiality of
	   this implicit game is \us{difficult since it requires access to
		   a} closed-form expression
	   for $ \Sscr(\cdot) $.  \us{Motivated by this, we} now provide a different
	   avenue that relies on the property
	   of {\em quasi-potentiality}.
\begin{definition} [Quasi-potential multi-leader multi-follower games] \label{def:quasipot}
	Consider a multi-leader multi-follower game $\Epec$ in which the
	player objectives are denoted by $\{\varphi_1, \hdots,
	\varphi_N\}$. 
	$\Epec$ is referred to as a quasi-potential
	game if the following hold:
	\begin{enumerate}
\item[(i)] For $i \in \Nscr$, there exist functions
	$\phi_1(x),\hdots,\phi_N(x)$ and a function $h(x,y_i)$ such that each player $i$'s objective $\varphi_i(\cdot)$ is given as
	$ \varphi_i(x_i,y_i;x\mi) \equiv \phi_i(x) + h(x,y_i).  $
	\item[(ii)] There exists a function $\pi(\cdot)$ such
	that for all $i = 1, \hdots, N$, and for all $x\in X$ and $x_i'\in X_i$, we have 
	$\phi_i(x_i;x^{-i}) - \phi_i(x_i';x^{-i}) = \pi(x_i;x^{-i}) -
	\pi(x'_i;x^{-i}).$ \\
	The function $\pi +h$ is called the quasi-potential function.
	\end{enumerate}
\end{definition}
Notice that the function $ h $ does not have a subscript `$ i $', and is thereby the same for all leaders. 
Essentially a quasi-potential game has objective functions which can be a written as a sum of an `$ x $' part and an `$x,y_i$' part, wherein the `$ x $' part admits a potential function in the standard sense and the `$ x,y_i $' part is identical for all players.
There exists a function $ \pi $ as required in Definition
\ref{def:quasipot} if and only if~\cite{monderer96potential}  for all $x \in \prod_{i=1}^N X_i$
\begin{equation}
 \nabla_{x_i} \pi(x) \equiv \nabla_{x_i}\phi_i(x) \qquad \forall i\in \Nscr. \label{eq:phi}
\end{equation}
\begin{remarkc} \label{rem:pot}
A potential multi-leader multi-follower  game is one where leaders have objective
functions $\varphi_i, i \in \Nscr$ such that there exists a function $\pi$, called potential function,
	 such that for all $i\in \Nscr$, for all $(x_i,x\mi) \in X,
	 (y_i,y\mi) \in Y$ and for all $x_i' \in X_i,y_i' \in Y_i,$
$\varphi_i(x_i,y_i;x\mi,y\mi)-\varphi_i(x_i',y_i';x\mi,y\mi) = \pi(x_i,y_i;x\mi,y\mi) - \pi(x_i',y_i';x\mi,y\mi). \label{eq:pot}$
Notice that a quasi-potential game is a potential game, but the converse is not true.
\end{remarkc}

Consider the  optimization problem 
   P$^{\rm quasi}$, defined as:
\begin{align}
&{	\problemcompactsmall{P$^{\rm quasi}$}
	{x,w}
	{\pi(x) + h(x,w)}
				 {(x,w) \in \Fscr^{\rm quasi},}} \non \\
				 & \where \non\\
&\Fscr^{\rm quasi} \triangleq \left\{(x,w) \ \left\lvert \
		 x_i \in X_i, \ i \in \Nscr, w \in \Sscr(x)
\right. \right\}. \label{eq:fscrcent} \end{align}
 
Our main result relates the minimizers of
	P$^{\rm quasi}$ to the global equilibria of $\Epec$. 
\begin{proposition}\label{unshared-relation}
	Consider a quasi-potential multi-leader multi-follower game $\Epec$. Then
	if $(x,w)$ is a global minimizer of {\rm P}$^{\rm quasi}$, then $(x,y)$, where $y_i=w$ for all $i \in \Nscr,$ is a global
	equilibrium of $\Epec$.
\end{proposition}
\begin{IEEEproof}
Observe that a point $(x_i,w)$ is feasible for the 
	$i\th$ agent's problem  L$_i(x^{-i})$ if and only if $ (x,w)\in \Fscr^{\rm quasi} $. \ie, for each $ i $ 
\begin{equation}
 (x_i,w) \in \Omega_i(x\mi)\quad \iff \quad(x,w)\in \Fscr^{\rm quasi}.\label{eq:potiff}
\end{equation}  
Now suppose, $(x,w)$ is a solution  of P$^{\rm quasi}$. Then, 
	$$ \pi(x) + h(x,w) \leq \pi(x') + h(x',w'), \forall x' \in X$$
 and $ \forall w' \in S(x'). $
	More specifically, taking $x'=(x_i',x\mi)$ and $w' \in S(x'_i,x\mi)$, and using \eqref{eq:potiff}, it
	follows that
	$$ \pi(x_i;x_{-i}) + h(x,w) \leq \pi(x_i';x_{-i}) + h((x'_i,x_{-i}),w'),$$
	$\forall (x_i',w') \in \Omega_i(x\mi)$. 
	Since $\Epec$ is a quasi-potential game, 
	$$\pi(x_i,x_{-i})-\pi(x_i',x_{-i}) = \phi_i(x_i;x_{-i}) -
	\phi_i(x_i',x_{-i}).$$ Let $y$ be such that $y_i=w$ for all $i\in \Nscr$. 
	As a result
	$ \varphi_i(x_i,y_i;x\mi) \leq \varphi_i(x_i',w';x\mi),\ \forall x_i',w' \in \Omega_i(x\mi),$
and hence $(x_i,y_i)$ is a solution L$_i(x^{-i})$. The
	result follows.
\end{IEEEproof}

Given this relationship between the minimizers of P$^{\rm quasi}$ and the
equilibria of $\Epec$, existence of equilibria is guaranteed by the
solvability of P$^{\rm quasi}$, as formalized by the next result.
\begin{theorem}[Existence of global equilibria of $\Epec$] \label{thm:unshared-pot}
Let $\Epec$ be a quasi-potential multi-leader multi-follower game. Suppose $\Fscr^{\rm quasi}$ is a nonempty set and
$\varphi_i$ is a continuous function for $i = 1, \hdots, N$. If the
minimizer of {\rm P}$^{\rm quasi}$ exists
(for example, if either $\pi$ is a coercive function over $\Fscr^{\rm quasi}$ or if
 $\Fscr^{\rm quasi}$ is compact ), then $\Epec$ admits  an equilibrium. 
\end{theorem}

Observe that $ \Fscr^{\rm quasi} $ is a closed set if $ K,G $ are continuous: the constraint $ w \in \Sscr(x) $ is equivalent to a nonlinear equation
$\fnat(w;x) =0,$ 
where $ \fnat(\cdot;x) $ is the natural map~\cite{facchinei03finiteI} of $ \VI(G(\cdot;x),K(x))$, given by
$\fnat(y;x) =y - \Pi_{K(x)}(y - G(y;x)), $
where $ \Pi_{Z}(z) $ is the projection of $ z $ on a set $ Z. $
 By the continuity of $ \fnat $ the zeros of $ \fnat $ form a closed set. Consequently, if $ K,G $ are continuous (which would hold, \eg, when the follower objectives are continuously differentiable and their constraints are independent of $ x $), then 
ascertaining the compactness of $ \Fscr^{\rm quasi} $ amounts to only
its boundedness. More generally, {\rm P}$^{\rm quasi}$ is a mathematical
program with equilibrium constraints and under coercivity of the objective or compactness
of the feasible region, a global
minimizer exists~\cite[Ch.~1]{luo96mathematical}.

Consider a special case of $ \Epec $, denoted $\Epec^{\rm ind}$,  where the $i^{\rm th}$ leader solves
the following problem wherein the objective is independent of $y_i$.
$$
{	\problemcompactsmall{L$^{\rm ind}_i(x^{-i})$}
	{x_i,y_i}
	{\varphi_i(x_i;x^{-i})}
				 {\begin{array}{r@{\ }c@{\ }l}
		x_i &\in& X_i, \ 		y_i \in \Sscr(x).
	\end{array}}}
	$$
	The following corollary captures the relationship between the global
	minimizers of P$^{\rm quasi}$ and the global equilibria of
	$\Epec^{\rm ind}$. 
\begin{corollary}
	Consider game $\Epec^{\rm ind}$ in which for each $ i $, $\varphi_i(x_i,y_i;x\mi) \equiv \varphi_i(x_i,x\mi)$ and the functions 
	$ \{\varphi_i\}_{i\in \Nscr} $ admit a potential function. Then this game is a quasi-potential multi-leader multi-follower game. Further, 
	if $(x,w)$ is a global minimizer of {\rm P}$^{\rm quasi}$, then $ (x,y) $ where $ y_i=w  $ for all $ i\in \Nscr $ is a global
	equilibrium of $\Epec^{\rm ind}$. If a solution exists to {\rm P}$^{\rm quasi}$, it is a global equilibrium of $ \Epec^{\rm ind}.$
\end{corollary}
The proof of this corollary follows from noting that $ \Epec^{\rm ind} $ is trivially a quasi-potential game (in Definition \ref{def:quasipot}, take $ h \equiv 0$ and $ \pi $ to be the given potential function of $ \{\varphi_i\}_{i \in \Nscr} $).

{At this juncture, it is worth differentiating the above existence
statements from more standard results presented
	in~\cite{sherali84multiple,su07analysis} where the follower
	equilibrium decisions are eliminated by leveraging the
	single-valuedness of the solution set of the follower equilibrium
	problem. In these approaches, the final claim rests on showing that the
	implicitly defined objective function (in the $ x $-space) is convex and continuous,
	properties that again require further assumptions. In comparison,
	we do not impose any such requirement.} {Finally, we believe that
		the class of {\em quasi-potential games} is not an artificial
			construct. For instance, the congestion control games
			arising in communication networks in \ref{sec:examples} lead to a
			quasi-potential multi-leader multi-follower game}. 

\begin{remarkc} One may ask the following question: if the objectives of
the leaders admit a potential function in
the $ x,y $ space, as in Remark \ref{rem:pot},  then does the resulting  game have an equilibrium?
The answer is no, as demonstrated by an example of Pang and Fukushima
\cite{pang05quasi}. We consider this example in
Section~\ref{sec:pot-examp}. This poses a challenge for generalizing our
results beyond the class of quasi-potential games.  \end{remarkc}

\subsubsection{Pessimistic formulation}
In this formulation, the $i\th$ leader solves the
following problem:
$${
	\minmaxproblemcompactsmall{\underbar{L}$_i(x\mi)$}
	{x_i}{y_i}
	{\varphi_i(x_i,y_i;x\mi)}
				 {\begin{array}{r@{\ }c@{\ }l}
				 x_i &\in& X_i,\ y_i \in \Sscr(x)			
	\end{array}}}
	$$
	Denote the resulting game by \underbar{$ \Epec $}; an equilibrium of \underbar{$ \Epec $} is defined as a Nash equilibrium of $ \{\underbar{\rm L}_i\}_{i\in \Nscr}. $ 
We know of no existence results for equilibria of the pessimistic formulation. Indeed, solutions to individual leader problems 
\underbar{L}$_i $ even may  not exist since each \underbar{L}$_i $ is
effectively a \textit{trilevel} optimization problem (the problem that
		defines $ \Sscr(\cdot),$ which is nested inside the maximization
		over $ y_i $, which in turn is nested inside the minimization
		over $ x_i $); for such problems, strong continuity properties
of inner nested problems are as good as
necessary~\cite{facchinei03finiteI} for the problem to admit a solution.
Of course, when $ \Sscr(\cdot) $ is single-valued, the pessimistic and
optimistic formulations coincide and thereby the Theorem
\ref{thm:unshared-pot} applies to \underbar{$ \Epec $}.  Nevertheless,
	if the leader objectives $ \varphi_i, i \in \Nscr $ admit a
	quasi-potential function, one can obtain a relation between the
	solution, it if exists, of a problem analogous to P$^{\rm quasi} $
	and the equilibirum of \underbar{$ \Epec $}. Consider the problem
$$
{	\minmaxproblemcompactsmall{\underbar{P}$^{\rm quasi}$}
	{x_i}{w}
	{\pi(x) + h(x,w)}
				 {\begin{array}{r@{\ }c@{\ }l}
				 x &\in& X,\ w \in \Sscr(x).			
	\end{array}}}
	$$
\begin{theorem}
Consider the game \underbar{$ \Epec $} and suppose that the objectives of the leaders admit a quasi-potential function. If 
$ (x,w) $ is a solution of \underbar{P}$^{\rm quasi} $ then $ (x,y) $ where $ y_i =w $ for all $ i\in \Nscr $ is an equilibrium of 
\underbar{$ \Epec $}. Consequently, if \underbar{P}$ ^{\rm quasi} $ admits a solution, \underbar{$ \Epec $} admits a solution.
\end{theorem}
\begin{IEEEproof}
Let $ (x,w) $ solve \underbar{P}$^{\rm quasi}.$ Therefore, clearly,
\[ \pi(x) + \max_{w \in \Sscr(x)} h(x,w) \leq \pi(x') + \max_{w' \in \Sscr(x')} h(x',w') \ \forall x' \in X.\]
Let $ i\in \Nscr $ and $ \bar{x}_i \in X_i $ be arbitrary and put $ x' = (\bar{x}_i;x\mi) $ in the inequality above to get
$$ \pi(x) + \max_{w \in \Sscr(x)} h(x,w) \leq \pi(\bar{x}_i,x\mi) + \max_{w' \in \Sscr(\bar{x}_i;x\mi)} h(\bar{x}_i,x\mi,w')$$  
$ \forall \bar{x}_i \in X_i.$
Since \underbar{$ \Epec $} is a quasi-potential game, we get
\begin{align*}
 \phi_i(x) &+ \max_{w \in \Sscr(x)}h(x,w) \leq \phi_i(\bar{x}_i;x\mi) \\ 
 &+\max_{w' \in \Sscr(\bar{x}_i;x\mi)}h(\bar{x}_i,x\mi,w') \ \forall \bar{x}_i \in X_i. 
\end{align*}
Since this holds for each $ i, $ it follows that $(x,y) $ where $ y_i =w \ \forall  i\in \Nscr $ is an equilibrium of \underbar{$ \Epec. $}
\end{IEEEproof}
\begin{remarkc}\underbar{P}$^{\rm quasi}$
	represents an instance of a Stackelberg equilibrium problem and its
		solvability has been studied in several places, \eg,~\cite{basar99dynamic}.
\end{remarkc}
\subsection{Local and Nash stationary equilibria}\label{sec:local_quasi}
While the discussion thus far provides an approach for claiming
existence of global equilibria by obtaining a global solution to a
suitable optimization problem. However, the computation of a global
minimizer of P$^{\rm quasi}$ is a difficult
nonconvex problem that falls within the category of mathematical
programs with equilibrium constraints (MPEC). However, one can often
obtain stationary points or local minimizers of such problems and in
this subsection, we relate these points to analogous local or
stationarity variants of Nash equilibria.  {We begin with a formal
	definition of a Nash Bouligand stationary or a Nash B-stationary
		point.}\footnote{A 	primal-dual characterization of
			B-stationarity is provided by Pang and
				Fukushima~\cite{pang98complementarity}.} 
	
\begin{definition}[Nash B-stationary point]
A point $(x,y) \in \Fscr$ is a {\em Nash B-stationary point} of $\Epec$ if for all $i \in \Nscr$,
$\nabla_i\varphi_i(x,y)\t d \geq 0, \ \forall d \in \Tscr((x_i,y_i);\Omega_i(x\mi)), $
where $\Tscr(z;K)$, the tangent cone at $z \in K \subseteq
		\mathbb{R}^n$, is defined as follows:
\begin{align*}
\Tscr(z;K) \triangleq &\left\{ dz \in \Real^n: \exists
		0<\{\tau_k\} \to 0, K\ni\{z_k\}\to z \right.\\  
&	\mbox{ such that } dz = \lim_{k \to \infty}
		\left(({z_k - z})/{\tau_k}\right)\}. 
\end{align*} 
\end{definition}

	\begin{proposition} [Nash B-stationary points of $\Epec$]
	Consider a quasi-potential \mlmfg $\Epec$ and suppose $\varphi_i$ is
	a continuously differentiable function over $X \times Y$ for $i =
	1, \hdots, N$. If $(x,w)$ is a B-stationary point of
	${\rm P}^{\rm quasi}$, then $(x,y)$ where $ y_i =w $ for all $ i\in \Nscr $ is a Nash  B-stationary point of $\Epec$.
	\label{prop:stat}
	\end{proposition}
		\begin{IEEEproof}
	 A stationary point $(x,w)$ of P$^{\rm quasi}$ satisfies 
	\begin{align} \label{stat-piquasi}
		\nabla_{x}  (\pi(x)+h(x,w))\t d{x}  + \nabla_{w} h(x,w)
			\t  d{w} \geq 0, 
		\end{align}
$\forall (d{x},d{w}) \in
					\Tscr ((x,w);\Fscr^{\rm quasi})$.
	Fix an $i \in \Nscr$ and consider an arbitrary $(dx_i',dy_i') \in
	\Tscr(x_i,w;\Omega_i(x^{-i}))$. By the definition of the tangent cone, there exists a sequence $\Omega(x^{-i}) \ni( u_{i,k},v_{i,k})\buildrel{k}\over\rightarrow (x_i,w)$ and a sequence $0<\tau_k \buildrel{k}\over\rightarrow0$ such that $\frac{u_{i,k} -x_i}{\tau_k} \buildrel{k}\over\rightarrow dx_i'$ and $ \frac{v_{i,k}-w}{\tau_k}\buildrel{k}\over\rightarrow dy_i'$. It follows that the sequence 
		$ (\textbf{x}_{i,k}, \textbf{y}_{i,k}),$
		where $\textbf{x}_{i,k} = (x_1, \hdots, u_{i,k}, \hdots, x_N), \aur  \textbf{y}_{i,k} = v_{i,k},$
satisfies $(\textbf{x}_{i,k},\textbf{y}_{i,k}) \in \Fscr^{\rm quasi}$. Therefore, the direction $(\textbf{dx}_i,\textbf{dw})$ where 
		$ \textbf{dx}_i = (0, \hdots, dx_i', \hdots, 0) \mbox{ and } \textbf{dy}_i =
		 dy_i', $
		belongs to $\Tscr((x,w);\Fscr^{\rm quasi})$. 
Substituting $(dx,dw) = (\textbf{dx}_i,\textbf{dy}_i)$ in \eqref{stat-piquasi} and using \eqref{eq:phi} gives
\[\nabla_{x_i}(\phi_i(x)+h(x,w))\t dx_i' + \nabla_{w} h(x,w)\t dy_i'\geq 0.\]
Since, $i \in \Nscr$ and $(dx_i',dy_i') \in
	\Tscr(x_i,w;\Omega_i(x^{-i}))$ were arbitrary, $(x,w)$ is a Nash B-stationary point of $\Epec$. 
		\end{IEEEproof}
We now define a local Nash equilibrium and show
	its relationship to the local minimum of P$^{\rm quasi}$. 
\begin{definition}[Local Nash equilibrium] \label{def:local}
A point $(x,y) \in \Fscr$ is a {\em local Nash equilibrium} of $\Epec$ if for all $i \in \Nscr$, $(x_i,y_i)$ is a local minimum of {\rm L}$_i(x\mi)$.
\end{definition}

\begin{proposition} [Local Nash equilibrium of $\Epec$]
	Consider a quasi-potential \mlmfg $\Epec$ and suppose $\varphi_i$ is 
	continuously differentiable function over $X \times Y$ for $i = 1,
	\hdots, N$. If
	$(x,w)$ is a local minimizer of
	${\rm P}^{\rm quasi}$, then $(x,y)$, where  $ y=(w,\hdots,w) $ is a local Nash equilibrium of $\Epec$.
	\label{prop:local}
	\end{proposition}
\begin{IEEEproof}
If $(x,w)$ is a local minimum of P$^{\rm quasi}$, there exists a neighborhood of $(x,w)$, denoted
by $\Bscr(x,w)$, such that 
\begin{equation}
\pi(x) + h(x,w) \leq \pi(x') +h(x',w'), \label{eq:locquasi}
\end{equation}
for all $(x',w') \in \Bscr(x,w) \cap \Fscr^{\rm quasi}.  $
Consider an arbitrary $i \in \Nscr$ and let $\Bscr_i(x_i,w;x\mi)
:= \left\{(u_i,v_i)\ |  (u_i,x\mi,v_i) \in \Bscr(x,w)\right\}$.  Then by the definition of $ \Fscr^{\rm quasi}$ in 
\eqref{eq:fscrcent} it follows that
\begin{align*}
(u_i,v_i) &\in \Omega_i(x\mi) \cap
		\Bscr_i(x_i,w;x\mi)\\
& \iff (u_i,x\mi,v_i) \in \Fscr^{\rm quasi} \cap \Bscr(x,w).
\end{align*}Thus,  in \eqref{eq:locquasi} put $ x'= (u_i,x\mi), w'=v_i$ to get 
$$ \pi(x) + h(x,w) \leq \pi(u_i,x\mi) + h(u_i,x\mi,v_i),$$
$\forall  (u_i,v_i) \in \Omega_i(x\mi) \cap
		\Bscr_i(x_i,w;x\mi)$
 Then employing 
\eqref{eq:pot}, we get 
$$\varphi_i(x,w) \leq \varphi_i(u_i,x\mi,v_i), \ \forall\
		(u_i,v_i) \in\Omega_i(x\mi) \cap
				\Bscr_i(x_i,w;x\mi).$$
In other words, $(x_i,w)$ is a local minimizer of L$_i(x\mi)$. This holds for each $i \in \Nscr,$ whereby $(x,y)$ where $ y=(w,\hdots,w) $ is a local Nash equilibrium. 
\end{IEEEproof}
For sake for brevity we have chosen to focus only on the Nash
B-stationary points and local Nash equilibria. Similar relationships
hold also between other notions of stationarity.
\subsection{Revisiting the example of Pang and Fukushima~\cite{pang05quasi}} \label{sec:pot-examp}
In 2005, Pang and Fukushima~\cite{pang05quasi} presented an example of a
simple multi-leader multi-follower game  that had no equilibrium.  This
game has two leaders with objectives $\varphi_1$ and $\varphi_2$,
	 defined as follows: 
$ \varphi_1(x_1,y_1) = \half x_1 + y_1 \ {\rm and } \ \varphi_2(x_2,y_2) = -\half
x_2 - y_2. $ 
This game admits  a potential function in the $ x,y $ space, but the game is not a quasi-potential game. 
However, several variants of this game with 
a modified $\varphi_2$ have an equilibrium, since they turn out to be quasi-potential games. 
For example, in a particular variant, the
leader objectives are given as follows: 
$ \varphi_1(x_1,y_1) = \half x_1 + y_1 \quad {\rm and } \quad \varphi_2(x_2,y_2) = -\half
x_2 + y_2, $
and this game has an equilibrium. 
We discuss the original example and the variants 
next.
 
\begin{examplec}[Pang and Fukushima~\cite{pang05quasi}]\label{ex:31} In this example, $ \Nscr = \{1,2\} $
$X_1 = X_2 = [0,1]$ and $Y = \Real$, and there is one follower~\cite{pang05quasi} which solves the optimization problem
$$ \min_{y \geq 0} \left\{y(-1+x_1+x_2) + \half y^2\right\} = \max
\left\{ 0, 1-x_1-x_2\right\} $$  
Thus the problems $ {\rm L}_1, {\rm L}_2$ are as follows.
\begin{align*}
{\problemcompactsmall{L$_1(x_2)$}
	{x_1,y_1}
	{\half x_1 +  y_1 }
				 {\begin{array}{r@{\ }c@{\ }l}
				 x_1 & \in & [0,1], \\
				 y_1 &=&  \max \left\{ 0, 1-x_1-x_2\right\}.	
				 \end{array} }} \\ 
{\problemcompactsmall{L$_2(x_1)$}
	{x_2,y_2}
	{-\half x_2 -y_2 }
				 {\begin{array}{r@{\ }c@{\ }l}
				 				 x_2 & \in & [0,1], \\
				 				 y_2 &=&  \max\left\{ 0, 1-x_1-x_2\right\}.	
				 				 \end{array} }} 
\end{align*}
It can be easily checked that this game has no equilibrium. We refer the reader to \cite{pang05quasi} for details.

Observe that this game admits a potential function in the $(x,y)$ space given readily by 
$\pi(x,y) = \varphi_1(x_1,y_1) + \varphi_2(x_2,y_2) = \half x_1 +  y_1 -\half x_2 -y_2.$
However this game is \textit{not} a {quasi-potential game}, since,
clearly, one cannot find a function $h$ such that $h(x,y_1) = y_1$ and
$h(x,y_2)=-y_2$  to meet the requirement of Definition
\ref{def:quasipot}.  We now consider a modification of this
example.

\paragraph{Quasi-potential variants of Pang and Fukushima~\cite{pang05quasi}}
Consider the following variant of the Pang and Fukushima example: take $ \varphi_1(x_1,y_1) = \half x_1 +  h(y_1),  $ 
$ \varphi_2(x_2,y_2) = -\half x_2 +h(y_2), $
where $h(\cdot)$ is a continuous function. This game is a quasi-potential game, with quasi-potential function given by $\pi(x) + h(y) = \half x_1 -\half x_2 +h(y).$
By Proposition \ref{unshared-relation}, a global minimizer of $\Pcent$
is a global equilibrium of the $\Epec$, where $\Pcent$ is defined as 
\[{\problemsmall{$\Pcent$\,}
	{x,w}
	{\half x_1 -\half x_2 +h(w) }
				 {\begin{array}{r@{\ }c@{\ }l}
w &=& \max\left\{ 0, 1-x_1-x_2\right\},	\\
 				 x_1,x_2 &\in&  [0,1].\end{array}} }\]
Let us consider some special cases of $h$.
 
Take $h(w) \equiv w.$ In this case, in L$_1$, L$_2$, one may
substitute $y_1$ and $y_2$, resulting in leader problems $\tilde{\rm L}_1,\tilde{\rm L}_2$ in
$x_1,x_2$, with objectives, 
$\tilde{\varphi}_1(x_1,x_2) = \half x_1 + \max\{0,1-x_1-x_2\},\ \tilde{\varphi}_2(x_1,x_2) = -\half x_1 + \max\{0,1-x_1-x_2\},$
respectively. Notice that $\tilde{\varphi}_1$ and  $\tilde{\varphi}_2$ are convex in $x_1$ and
$x_2$ respectively. Since the feasible regions of both problems \us{are}
convex and compact (they are unit intervals) and the objective is convex
and continuous, \us{classical results
suffice for claiming~\cite{basar99dynamic}} that this game has an equilibrium.

 Take $h(w)\equiv-w$. In this case, again one may substitute for $y_1,y_2$ in terms of $x_1,x_2$. The resulting problems $\tilde{\rm L}_1$ and $\tilde{\rm L}_2$ are nonconvex.  Notice that $\Pcent$
is equivalent to minimizing $\half x_1 -\half x_2 - \max \{0,
   1-x_1-x_2\}$ over 
$\{(x_1,x_2): (x_1,x_2) \in [0,1]^2\}.$ It can be observed
   that the minimizer
of $\Pcent$ is given by $(x_1,x_2,w)=(0,0,1)$.
To see why the point
   $(x_1,x_2,y_1,y_2)=(0,0,1,1)$ is an equilibrium, notice that
   given $x_2=0,y_2=1$, the global minimizer of L$_1(x_2,y_2)$ is $x_1=0, y_1=1$. Similarly,
   with $x_1=0,y_1=1$, the global minimizer of L$_2(x_1,y_1)$ is again given by $x_2=0,y_2=1.$
\end{examplec}

Thus we see that quasi-potentiality is a powerful property that allows one to leverage the unique structure of multi-leader multi-follower games to claim existence of equilibria.
\section{Conclusions}\label{conc}
We consider multi-leader multi-follower games and
examine the question of the existence of an equilibrium. A standard
approach requires ascertaining when the reaction map admits fixed
points. However, this avenue  has several hindrances, an important one
being  the lack of continuity in the solution set associated with  the
equilibrium constraints capturing the follower equilibrium. {We observed
	that these challenges can be circumvented for {\em quasi-potential}
	multi-leader multi-follower games.
We show that any global minimizer of a suitably defined optimization
	problem is a global equilibrium of the game and that a similar
	result holds for the pessimistic formulation.  
Consequently, the above results reduced a question of the existence of
an equilibrium to that of the solvability of an optimization problem,
which can be claimed \us{under fairly standard conditions that are
	tractable and verifiable} -- \eg, coercive objective over a nonempty
	feasible region -- and the existence of a global equilibrium was
	seen to follow.  We further showed that local minima and
	B-stationary points of the respective MPECs are local Nash
	equilibria and Nash B-stationary points of the corresponding
	multi-leader multi-follower game.
%

\bibliographystyle{../Mylatexfiles/plainini}
\bibliography{../Mylatexfiles/ref,ref2}
\end{document}